\documentclass[11pt]{amsart}
\usepackage{amsopn}
\usepackage{amssymb, amscd}
\usepackage{multirow}
\usepackage{graphicx, graphics, epsfig}
\usepackage{faktor} 
\usepackage{enumerate}
\usepackage{tikz}
\usepackage{pgfplots}
\usepackage{hyperref}
\usepackage{color}
\usepackage{dynkin-diagrams}
\usepackage{epigraph}
\newcommand{\nc}{\newcommand}

\nc{\fg}{\mathfrak{f} } \nc{\vg}{\mathfrak{v} } \nc{\wg}{\mathfrak{w} }
\nc{\zg}{\mathfrak{z} } \nc{\ngo}{\mathfrak{n} } \nc{\kg}{\mathfrak{k} }
\nc{\mg}{\mathfrak{m} } \nc{\bg}{\mathfrak{b} } \nc{\ggo}{\mathfrak{g} } \nc{\eg}{\mathfrak{e} }
\nc{\ggob}{\overline{\mathfrak{g}} } \nc{\sog}{\mathfrak{so} }
\nc{\sug}{\mathfrak{su} } \nc{\spg}{\mathfrak{sp} } \nc{\slg}{\mathfrak{sl} }
\nc{\glg}{\mathfrak{gl} } \nc{\cg}{\mathfrak{c} } \nc{\rg}{\mathfrak{r} }
\nc{\hg}{\mathfrak{h} } \nc{\tg}{\mathfrak{t} } \nc{\ug}{\mathfrak{u} }
\nc{\dg}{\mathfrak{d} } \nc{\ag}{\mathfrak{a} } \nc{\pg}{\mathfrak{p} }
\nc{\sg}{\mathfrak{s} } \nc{\affg}{\mathfrak{aff} } \nc{\qg}{\mathfrak{q} } \nc{\lgo}{\mathfrak{l} } \nc{\holg}{\mathfrak{hol} } \nc{\stabg}{\mathfrak{stab} } \nc{\iautg}{\mathfrak{iaut} } \nc{\bautg}{\mathfrak{baut} } \nc{\hiautg}{\mathfrak{hiaut} }

\nc{\pca}{\mathcal{P}} \nc{\nca}{\mathcal{N}} \nc{\lca}{\mathcal{L}}
\nc{\oca}{\mathcal{O}} \nc{\mca}{\mathcal{M}} \nc{\tca}{\mathcal{T}}
\nc{\aca}{\mathcal{A}} \nc{\cca}{\mathcal{C}} \nc{\gca}{\mathcal{G}}
\nc{\sca}{\mathcal{S}} \nc{\hca}{\mathcal{H}} \nc{\bca}{\mathcal{B}}
\nc{\dca}{\mathcal{D}} \nc{\eca}{\mathcal{E}} \nc{\wca}{\mathcal{W}} \nc{\ica}{\mathcal{I}} \nc{\kca}{\mathcal{K}} \nc{\acca}{\mathcal{AC}}

\nc{\vp}{\varphi} \nc{\ddt}{\tfrac{d}{dt}} \nc{\dsdt}{\tfrac{d^2}{dt^2}} \nc{\dds}{\frac{d}{ds}}
\nc{\dpar}{\frac{\partial}{\partial t}} \nc{\im}{\mathrm{i}} \nc{\dpars}{\tfrac{\partial}{\partial t}} 

\nc{\SO}{\mathrm{SO}} \nc{\Spe}{\mathrm{Sp}} \nc{\Sl}{\mathrm{SL}}
\nc{\SU}{\mathrm{SU}} \nc{\Or}{\mathrm{O}} \nc{\U}{\mathrm{U}} \nc{\Gl}{\mathrm{GL}}
\nc{\Se}{\mathrm{S}} \nc{\Cl}{\mathrm{Cl}} \nc{\Spin}{\mathrm{Spin}}
\nc{\Pin}{\mathrm{Pin}} \nc{\G}{\mathrm{GL}_n(\RR)} \nc{\g}{\mathfrak{gl}_n(\RR)}
\nc{\Eg}{\mathrm{E}} \nc{\Fg}{\mathrm{F}} \nc{\Gg}{\mathrm{G}} \nc{\Gr}{\mathrm{Grass}}

\nc{\RR}{{\Bbb R}} \nc{\HH}{{\Bbb H}} \nc{\CC}{{\Bbb C}} \nc{\ZZ}{{\Bbb Z}}
\nc{\FF}{{\Bbb F}} \nc{\NN}{{\Bbb N}} \nc{\QQ}{{\Bbb Q}} \nc{\PP}{{\Bbb P}} \nc{\OO}{{\Bbb O}}

\nc{\vs}{\vspace{.3cm}} \nc{\vsp}{\vspace{1cm}} \nc{\ip}{\langle\cdot,\cdot\rangle}
\nc{\ipp}{(\cdot,\cdot)} \nc{\la}{\langle} \nc{\ra}{\rangle} \nc{\unm}{\tfrac{1}{2}}
\nc{\unc}{\tfrac{1}{4}} \nc{\und}{\frac{1}{16}} \nc{\no}{\vs\noindent}
\nc{\lam}{\rho^2(\RR^n)^*\otimes\RR^n} \nc{\tangz}{{\rm T}^{\rm Zar}}
\nc{\nor}{{\sf n}}  \nc{\mum}{/\!\!/} \nc{\kir}{/\!\!/\!\!/}
\nc{\Ri}{\tfrac{4\Ric_{\mu}}{||\mu||^2}} \nc{\ds}{\displaystyle}
\nc{\ben}{\begin{enumerate}} \nc{\een}{\end{enumerate}} \nc{\f}{\frac}
\nc{\lb}{[\cdot,\cdot]} \nc{\isn}{\tfrac{1}{||v||^2}}
\nc{\gkp}{(\ggo=\kg\oplus\pg,\ip)} \nc{\ukh}{(\ug=\kg\oplus\hg,\ip)}
\nc{\tgkp}{(\tilde{\ggo}=\kg\oplus\pg,\ip)}
\nc{\wt}{\widetilde}
\nc{\iop}{\mathtt{i}} \nc{\jop}{\mathtt{j}} 
\nc{\Hk}{H_{\kil}} \nc{\gk}{g_{\kil}}

\nc{\Hess}{\operatorname{Hess}} \nc{\ad}{\operatorname{ad}}
\nc{\Ad}{\operatorname{Ad}} \nc{\rank}{\operatorname{rk}}
\nc{\Irr}{\operatorname{Irr}} \nc{\End}{\operatorname{End}}
\nc{\Aut}{\operatorname{Aut}} \nc{\Inn}{\operatorname{Inn}}
\nc{\Der}{\operatorname{Der}} \nc{\Ker}{\operatorname{Ker}}
\nc{\Iso}{\operatorname{Isom}} \nc{\Diff}{\operatorname{Diff}}
\nc{\Lie}{\operatorname{L}} \nc{\tr}{\operatorname{tr}} \nc{\dif}{\operatorname{d}}
\nc{\sen}{\operatorname{sen}} \nc{\modu}{\operatorname{mod}}
\nc{\CRic}{\operatorname{PP}} \nc{\Cric}{\operatorname{P}} \nc{\Ricci}{\operatorname{Ric}}
\nc{\sym}{\operatorname{sym}} \nc{\herm}{\operatorname{herm}} \nc{\symac}{\operatorname{sym^{ac}}}
\nc{\symc}{\operatorname{sym^{c}}} \nc{\scalar}{\operatorname{Sc}}
\nc{\grad}{\operatorname{grad}} \nc{\ricci}{\operatorname{Rc}} \nc{\kil}{\operatorname{B}} \nc{\cas}{\operatorname{C}} \nc{\lic}{\operatorname{L}}
\nc{\Nor}{\operatorname{Norm}}  \nc{\ricc}{\operatorname{Rc^{c}}}
\nc{\Ricc}{\operatorname{Ric^{c}}} \nc{\ricac}{\operatorname{Rc^{ac}}}
\nc{\Ricac}{\operatorname{Ric^{ac}}} \nc{\Riem}{\operatorname{Rm}} \nc{\Sec}{\operatorname{Sec}}
\nc{\riccig}{\operatorname{ric^{\gamma}}} \nc{\mm}{\operatorname{m}} \nc{\Mm}{\operatorname{M}}
\nc{\Le}{\operatorname{L}} \nc{\tang}{\operatorname{T}}
\nc{\level}{\operatorname{level}} \nc{\rad}{\operatorname{r}}
\nc{\abel}{\operatorname{ab}} \nc{\CH}{\operatorname{CH}} \nc{\Cone}{{\mathcal C}} \nc{\CCone}{\operatorname{CC}} \nc{\CP}{{\mathcal P}}
\nc{\mcc}{\operatorname{mcc}} \nc{\Adj}{\operatorname{Adj}}
\nc{\Order}{\operatorname{O}}  \nc{\inj}{\operatorname{inj}} \nc{\proy}{\operatorname{pr}}
\nc{\vol}{\operatorname{vol}} \nc{\Diag}{\operatorname{Dg}} \nc{\Diagg}{\operatorname{Diag}}
\nc{\Spec}{\operatorname{Spec}} \nc{\Ima}{\operatorname{Im}} \nc{\Rea}{\operatorname{Re}}
\nc{\spann}{\operatorname{span}} \nc{\Aff}{\operatorname{Aff}} \nc{\E}{\operatorname{E}} \nc{\id}{\operatorname{id}} \nc{\dete}{\operatorname{det}} \nc{\Crit}{\operatorname{Crit}} \nc{\val}{\operatorname{val}} \nc{\Bihol}{\operatorname{Bihol}} 

\theoremstyle{plain}
\newtheorem{theorem}{Theorem}[section]

\theoremstyle{definition}

\theoremstyle{remark}
\newtheorem{remark}[theorem]{Remark}

\title{Low dimensional C-spaces}

\author{Jorge Lauret} \author{Facundo Montedoro} 

\address{FaMAF, Universidad Nacional de C\'ordoba and CIEM, CONICET (Argentina)} 
\email{jorgelauret@unc.edu.ar} 
\address{Department of Mathematics, Aarhus University (Denmark)}
\email{montedoro@math.au.dk} 

\thanks{This research was partially supported by a grant from Universidad Nacional de C\'ordoba (Argentina).  }

\date{\today}

\begin{document}

\maketitle

\begin{abstract} 
We give the classification of C-spaces (i.e., compact homogeneous spaces admitting an invariant complex structure), up to finite cover and equivariant diffeomorphism, of dimension $6$, $8$ and $10$.  The conditions for the existence of a balanced, pluriclosed, CYT and LCK Hermitian metric are established for dimensions $6$ and $8$.    
\end{abstract}

\tableofcontents


\section{Introduction}\label{intro}

A compact homogeneous space $M=G/K$ admitting a $G$-invariant complex structure is called a {\it C-space}.  We know from results by Wang \cite{Wng}, Samelson \cite{Sml} and Tits \cite{Tts} (see also \cite{NiWll}), that any C-space can be constructed according to the following recipe, whose only ingredients are a flag manifold $F=G_f/H_f$ and a compactly embedded subspace $\zg(\kg)\subset\zg(\hg_f)$, where $\zg(\hg_f)$ is the center of the Lie algebra $\hg_f$ of $H_f$:

\begin{quote}
$G=G_f\times T^z$ and the C-space is given by
$$
M=G/K=G_f/K\times T^z,
$$ 
where $K\subset H_f$ is the Lie subgroup with Lie algebra $[\hg_f,\hg_f]\oplus\zg(\kg)$.  
\end{quote}
Note that $M$ is therefore a torus fibration over $F$ (see \eqref{fib}).  
 
The reductive decomposition $\ggo_f=\hg_f\oplus\qg$ for the flag $F=G_f/H_f$ determines the following reductive decomposition for the C-space $M=G/K$: 
$$
\ggo= \kg\oplus\pg = 
\underbrace{[\hg_f,\hg_f]\oplus\rlap{$\overbrace{\phantom{\zg(\kg) \oplus\ag}}^{\zg(\hg_f)}$}\zg(\kg)}_\kg \oplus \underbrace{\ag_f\oplus\zg(\ggo)\oplus\qg}_\pg,  \qquad \ag:=\ag_f\oplus\zg(\ggo),  
$$
where everything is $Q$-orthogonal relative to a fixed bi-invariant inner product $Q$ on $\ggo$.  

It is proved in \cite{NiWll} that any $G$-invariant complex structure on $M=G/K$ descending to the flag $F$ is of the form $J=(J_\ag,J_\qg)$, 
where $J_\ag:\ag\rightarrow\ag$ is any linear map such that $J_\ag^2=-I$ and $J_\qg:\qg\rightarrow\qg$ is the complex structure on $F$ attached to some choice $\Delta_\qg^+$ of positive complementary roots of $F$; namely, $J_\qg E_\alpha = \pm\im E_\alpha$ for all $E_\alpha\in\ggo_\alpha$, $\alpha\in\pm\Delta_\qg^+$.  

In this paper, we obtain the following classification results.  

\begin{theorem}\label{clasif}
The lists of C-spaces, up to finite cover and equivariant diffeomorphism, of dimension $6$, $8$ and $10$ are as in Tables \ref{tab6}, \ref{tab8} and \ref{tab10}, respectively.  
\end{theorem}

See \S\ref{cla-sec} and also Tables \ref{tab8-2} and \ref{tab10-2} for further information on these C-spaces.  A partial classification of $12$-dimensional C-spaces is given in Table \ref{tab12}.  

In \S\ref{exi-sec}, we study the Hermitian geometry of invariant complex structures.  

\begin{theorem}\label{clasif2}
The $G$-invariant complex structures on each $6$- and $8$-dimensional C-space $M=G/K$ descending to the corresponding flag $F$ are respectively given in the third column of Tables \ref{tab6-exi} and \ref{tab8-exi}.  Furthermore, for each of these complex structures, the tables include the conditions for the existence of a $G$-invariant Hermitian metric of the following special kinds: balanced, pluriclosed (SKT), Calabi-Yau with torsion (CYT) and locally conformally K\"ahler (LCK).     
\end{theorem}

By way of example, we highlight the landscape of complex structures on the following C-spaces:  

\begin{enumerate}[{\small $\bullet$}]
\item The space of all left-invariant complex structures on the Lie group 
$$
M^8=\SU(2)\times\SU(2)\times T^2 
$$
consists, up to biholomorphism, of a single component $\cca$, which is identified with the $8$-dimensional homogeneous space $\Gl_4(\RR)/\Gl_2(\CC)$ (see Tables \ref{tab8}, \ref{tab8-2} and \ref{tab8-exi}, line 2).  We have that 
$$
\cca_{bf}=\cca_{plu}\subsetneq\cca_{cyt}\subsetneq\cca, \qquad \dim{\cca_{plu}}= 2, \quad \dim{\cca_{cyt}}=8,
$$ 
where $\cca_{bf}$, $\cca_{plu}$ and $\cca_{cyt}$ are the subsets of $\cca$ of complex structures respectively admitting a Bismut-flat, pluriclosed and CYT compatible metric.  The space $\cca_{cyt}$ is the disjoint union of an open subset of $\cca$ and a $4$-dimensional subspace contained in its boundary, which is in turn open in the set of all product complex structures.        

\item The space $\cca$ of all invariant complex structures on the C-space 
$$
M^8=\SU(3)/S^1_{p,q}\times S^1, 
$$ 
where ${\rm AW}_{p,q}=\SU(3)/S^1_{p,q}$ is an Aloff-Wallach space, consists of three components: $\cca=\cca_1\sqcup\cca_2\sqcup\cca_3$ (see Tables \ref{tab8}, \ref{tab8-2} and \ref{tab8-exi}, line 8).  The existence of special compatible metrics on each component is described in Table \ref{tab-in}.  

\item For the C-space 
$$
M^8=S^3\times S^3\times S^2=\SU(2)^3/S^1_{p,q,r},
$$ 
given in Tables \ref{tab8}, \ref{tab8-2} and \ref{tab8-exi}, line 12, we have four components: $\cca=\cca_1\sqcup\cca_2\sqcup\cca_3\sqcup\cca_4$, and the existence of special compatible metrics is described in Table \ref{tab-in2}.  The subset $\cca_{cyt}$ is nonempty, proper and open in the corresponding component $\cca_i$.  By the symbol ? we mean that none complex structure in the component admits an special metric of the kinds studied in this paper.  
\end{enumerate}

\begin{table}
\centering
\begin{tabular}{c | c | c | c }
\hline\hline \rule{0pt}{14pt}
$(p,q)$ & $\cca_1$ & $\cca_2$ & $\cca_3$ \\ [2pt]
\hline\hline\rule{0pt}{14pt}
(1,1) & Bal & CYT and LCK & Bal  \\ 
\hline \rule{0pt}{14pt}
(1,0) & Bal & LCB  & Bal  \\ 
\hline \rule{0pt}{14pt}
$p>q>0$ & Bal & LCK & Bal  \\ 
\hline\hline 
\end{tabular}
\vspace{.3cm}
\caption{Special Hermitian metrics on the $8$-dimensional C-spaces $\SU(3)/S^1_{p,q}\times S^1$.  LCB: locally conformally balanced.}\label{tab-in}
\end{table}

\begin{table}
\centering
\begin{tabular}{c | c | c | c | c}
\hline\hline \rule{0pt}{14pt}
$(p,q,r)$ & $\cca_1$ & $\cca_2$ & $\cca_3$ & $\cca_4$ \\ [2pt]
\hline\hline\rule{0pt}{14pt}
$(1,1,0)$ & ? & $\cca_{cyt}$ & ? & $\cca_{cyt}$ \\ 
\hline \rule{0pt}{14pt}
$p=q+r$, $r>0$ & Bal & ? & ? & $\cca_{cyt}$ \\  
\hline \rule{0pt}{14pt}
$p\ne q+r$ & ? & ? & ? & ? \\ 
\hline\hline 
\end{tabular}
\vspace{.3cm}
\caption{Special Hermitian metrics on the $8$-dimensional C-spaces $\SU(2)^3/S^1_{p,q,r}$.}  
\label{tab-in2}
\end{table}

\section{Preliminaries}\label{preli}

In this section, we give an overview on C-spaces, their complex structures, their Hermitian metrics and the conditions for the existence of special metrics.  We refer to \cite{Herm} for a more detailed treatment.  

\subsection{C-spaces}\label{Csp-sec}
A compact homogeneous space $M=G/K$ is called a {\it C-space} if there exists a Lie subalgebra $\hg\subset\ggo$ such that $\hg$ is the centralizer of some abelian subalgebra of $\ggo$, $\kg\subset\hg$ and $[\hg,\hg]=[\kg,\kg]$, 
where $\ggo$ and $\kg$ are the Lie algebras of $G$ and $K$, respectively.  Each of these subalgebras $\hg\subset\ggo$ therefore determines a flag manifold $F=G/H$ and the following torus fibration over $F$, called a {\it Tits fibration}:  
\begin{equation}\label{fib}
A=H/K\longrightarrow M=G/K \longrightarrow F=G/H, 
\end{equation}
with fiber the torus $A=H/K$.  Given any background bi-invariant inner product $Q$ on $\ggo$, the $Q$-orthogonal reductive decomposition for a C-space $M=G/K$ is given by  
\begin{equation}\label{reddec}
\ggo= \kg\oplus\pg = 
\underbrace{[\hg,\hg]\oplus\rlap{$\overbrace{\phantom{\zg(\kg) \oplus\ag}}^{\zg(\hg)}$}\zg(\kg)}_\kg \oplus \underbrace{\ag\oplus\qg}_\pg,  \qquad \pg=\ag\oplus\qg\equiv T_oM, \quad T_oA\equiv\ag, \quad T_oF\equiv\qg, 
\end{equation}
where $[\ggo,\ggo]= (\hg\cap[\ggo,\ggo])\oplus\qg$ is the reductive decomposition of the flag $F$.  We can assume that $\zg(\ggo)\subset\ag$, thus the C-space $M=G/K$ is determined by only the flag $F=G_f/H_f$ and a compactly embedded subspace $\zg(\kg)\subset\zg(\hg_f)$, where $G_f:=[G,G]$ and $H_f:=H\cap[G,G]$.  In particular, $G_f/K$ is also a C-space and 
\begin{equation}\label{gprime}
M=G/K=G_f/K\times T^z, \qquad z:=\dim{\zg(\ggo)}.
\end{equation}   
In that case, $\ag=\zg(\ggo)\oplus\ag_f$, where $\ag_f:=\ag\cap[\ggo,\ggo]$, so $\zg(\hg_f)=\zg(\kg)\oplus\ag_f$.  

As a manifold, any C-space can also be decomposed as
\begin{equation}\label{dec}
M=G/K=\widetilde{G}/\widetilde{H}\times \underline{G}/\underline{K}\times\overline{G}, 
\end{equation}
where $\widetilde{G}/\widetilde{H}$ is a flag manifold, $\underline{G}/\underline{K}= [\underline{G},\underline{G}]/\underline{K}\times T^z$ is a {\it strict} C-space (see \cite[Definition 3.8]{Herm}) and $\overline{G}$ is a compact semisimple Lie group.   

As usual, we consider the root system $\Delta\subset(\tg^c)^*$ defined by the Cartan subalgebra $\tg^c:=\tg\otimes\CC$ of the complex semisimple Lie algebra $\ggo_f^c$, where $T\subset H_f\subset G_f$ is a maximal torus and $\ggo_f:=[\ggo,\ggo]$.  The $Q$-orthogonal reductive decomposition 
$\ggo_f=\hg_f\oplus\qg$ of $F=G_f/H_f$, is given by 
\begin{equation}\label{hq-intro}
\hg_f^c=\tg^c\oplus\bigoplus_{\alpha\in\Delta_{\hg_f}}\ggo_\alpha, \qquad 
\qg^c=\bigoplus_{\alpha\in\Delta_\qg}\ggo_\alpha, \qquad \Delta=\Delta_{\hg_f}\sqcup\Delta_\qg,
\end{equation}
where $\ggo_\alpha:=\{ E\in\ggo_f^c:[H,E]=\alpha(H)E,\; H\in\tg^c\}$, $\Delta_{\hg_f}:=\{\alpha\in\Delta:\alpha|_{\zg(\hg_f)}=0\}$ is the root system of $[\hg_f,\hg_f]$ and $\Delta_\qg$ is the set of {\it complementary roots}.

\subsection{Complex structures}\label{comp-sec}
Given an even dimensional C-space $M=G/K$ and one of its Tits fibrations as in \eqref{fib}, attached to a subalgebra $\hg\subset\ggo$, every $G$-invariant complex structure on $M=G/K$ descending to the flag $F$ is of the form 
\begin{equation}\label{J}
J=(J_\ag,J_\qg)=\left[\begin{matrix} J_\ag&0\\0&J_\qg\end{matrix}\right], 
\end{equation}
where $J_\ag:\ag\rightarrow\ag$ is any linear map such that $J_\ag^2=-I$ (in particular, $\dim{\ag}$ is even) and $J_\qg:\qg\rightarrow\qg$ is the complex structure on $F$ attached to some {\it invariant ordering} $\Delta_\qg^+$ of $\Delta_\qg$: 
$$
\Delta_\qg=\Delta_\qg^+\sqcup -\Delta_\qg^+, \quad (\Delta_{\hg_f}+\Delta_\qg^+)\cap\Delta\subset\Delta_\qg^+, \quad  (\Delta_\qg^++\Delta_\qg^+)\cap\Delta\subset\Delta_\qg^+, 
$$
that is, $J_\qg E_\alpha = \pm\im E_\alpha$ for all $E_\alpha\in\ggo_\alpha$, $\alpha\in\pm\Delta_\qg^+$.  The Tits fibration therefore becomes a holomorphic fibration: 
\begin{equation}\label{holfib}
(A,J_\ag)\longrightarrow (M=G/K,J) \longrightarrow (F,J_\qg),     
\end{equation}
and according to \cite[Theorem 1.3]{NiWll}, every $G$-invariant complex structure on $M=G/K$ is of the form $J=(J_\ag,J_\qg)$ for some of its Tits fibrations.  Moreover, up to biholomorphism, it is enough to consider only one of its Tits fibrations, which will be fixed from now on.  

The {\it Koszul vector} of $J$, defined by
\begin{equation}\label{koszul2}
Z_{J_\qg}:= \sum_{\alpha\in\Delta_\qg^+} \im H_\alpha\in\zg(\hg)=\zg(\kg)\oplus\ag, 
\end{equation}
where $\alpha(X)=Q(X,H_\alpha)$ for all $X\in\tg^c$, satisfies that the first Chern class $c_1(M,J)\in H^2(M)$ vanishes if and only if 
\begin{equation}\label{c10} 
Z_{J_\qg}\in\ag_f. 
\end{equation}
The set of all possible $J_\qg$ is finite, it is given by $\{J_1,\dots,J_k\}$, where $k$ is the number of chambers in 
\begin{equation}\label{ch}
\im\zg(\hg)\setminus\bigcup_{\alpha\in\Delta_\qg}\Ker\alpha\cap\im\zg(\hg), 
\end{equation}
or equivalently, the number of invariant orderings.  We thus obtain that the space $\cca$ of all $G$-invariant complex structures on $M=G/K$ descending to the flag $F$ is given by    
\begin{equation}\label{Ccal}
\cca=\cca_1\sqcup\dots\sqcup\cca_k, \qquad \dim{\cca_i}=2d^2,  \qquad \dim{\ag}=2d,
\end{equation}
where 
$$
\cca_i:=\left\{ (J_\ag,J_i): J_\ag:\ag\rightarrow\ag \;\mbox{linear},\;  J_\ag^2=-I\right\} \equiv \Gl_{2d}(\RR)/\Gl_d(\CC).
$$ 
Note that $c_1(J)=0$ for $J=(J_\ag,J_i)$ if and only if $Z_{J_i}\in\ag_f$ (see \eqref{c10}).

\subsection{Hermitian metrics}\label{herm-sec}
The isotropy representation $\pg$ of a C-space $M=G/K=G_f/K\times T^z$ as in \eqref{gprime} decomposes in $\Ad(K)$-invariant subspaces as
\begin{equation}\label{isotrep}
\pg=\ag\oplus\mg_1\oplus\dots\oplus\mg_r, \qquad \ag=\zg(\ggo)\oplus\ag_f, \qquad \qg=\mg_1\oplus\dots\oplus\mg_r, 
\end{equation}
where 
$$
\mg_j^c:=\bigoplus_{\alpha\in \Delta_j} \ggo_{\alpha}, \qquad \Delta_j:=\left\{\alpha\in \Delta:\alpha|_{\zg(\hg_f)}=\pm \alpha_j|_{\zg(\hg_f)}\right\}, 
$$ 
and $\{\alpha_1,\dots,\alpha_r\}\subset\Delta_\qg$ is a choice of {\it central roots} for the flag $F=G_f/H_f$, i.e.,  
\begin{equation}\label{cenrts}
\left\{\alpha|_{\zg(\hg)}:\alpha\in \Delta\right\}=\left\{\pm\alpha_1|_{\zg(\hg)},\dots,\pm\alpha_r|_{\zg(\hg)}\right\}.  
\end{equation}
A formula for the number $r$ is given in \cite{Fr}.  The action of $K$ on $\ag$ is trivial and the subspaces $\mg_j$ are always $\Ad(H_f)$-invariant, $\Ad(H_f)$-irreducible and pairwise $\Ad(H_f)$-inequivalent.  If the $\mg_j$ are also $\Ad(K)$-irreducible and pairwise $\Ad(K)$-inequivalent, which holds generically (see \cite[Remark 3.25]{Herm}), then any $G$-invariant metric on $M=G/K$ is given by 
\begin{equation}\label{g}
g=g_\ag+g_\qg, \qquad g_\qg=(x_\alpha)_{\alpha\in\Delta_\qg}=x_1Q|_{\mg_1}+\dots+x_rQ|_{\mg_r}, \quad x_i>0,
\end{equation}
where $g_\ag$ is any inner product on $\ag$.  Note that $(J,g)$ is compatible if and only if $(J_\ag,g_\ag)$ is compatible, and in that case, the corresponding K\"ahler form $\omega=g(J\cdot,\cdot)$ of the Hermitian structure $(J,g)$ is given by  
\begin{equation}\label{o}
\omega=\omega_\ag+\omega_\qg, \qquad \omega_\qg(E_\alpha,E_{-\alpha})=-\im x_\alpha, \quad\forall \alpha\in\Delta_\qg^+,   
\end{equation} 
and zero otherwise, where $\omega_\ag=g_\ag(J_\ag\cdot,\cdot)$.  The {\it metric Koszul vector}  
\begin{equation}\label{koszul-g}
Z_{J_\qg,g_\qg}:=\sum_{\alpha\in\Delta_\qg^+}\tfrac{1}{x_\alpha}  \im H_\alpha   
\in \zg(\hg_f)=\zg(\kg)\oplus\ag_f,  
\end{equation}
plays an important role in the geometry of the Hermitian C-space $(M=G/K,J,g)$.

\subsection{Special metrics}\label{sp-sec}
The following kinds of Hermitian metrics are considered distinguished on a given complex manifold:  

\begin{enumerate}[{\small $\bullet$}] 
\item {\it Balanced}: $d_g^*\omega=0$, or equivalently, $d\omega^{n-1}$=0 (see \cite{Pds, FinGrnVzz, Fin, FinGrn2, Kwn}).  

\item {\it Locally conformally balanced} (LCB): $d\theta_L=0$, where $\theta_L=-d_g^*\omega(J\cdot)$ is the Lee form (see \cite{FinTms, Prd}).  

\item {\it Pluriclosed}, or {\it Strong K\"ahler with torsion} (SKT): $\partial\overline{\partial}\omega=0$, or equivalently, $dd^c\omega=0$ (see \cite{Pds, FinGrnVzz, FinGrn1, FinGrn2, Fin, Brb, SKT-LG, Kwn}).  

\item {\it Calabi-Yau with torsion} (CYT): $\rho^b=0$, where $\rho^b$ is the Bismut Ricci form (see \cite{Grn, Pds, FinGrn1, FinGrn2}).  

\item {\it Locally conformally K\"ahler} (LCK): $d\omega=\theta\wedge\omega$ for some closed $1$-form $\theta$ (see \cite{OrnVrb, HsgKms1,HsgKms2, AlkCrtHsgKms, Gn, AlkHsgKms}).  If in addition $\nabla^g\theta=0$, the structure is called {\it Vaisman} (see \cite{Ist}).  
\end{enumerate}

Given a complex C-space $(M=G/K,J)$ without a flag factor (see \eqref{dec}), the conditions on $J=(J_\ag,J_\qg)$ for the existence of a compatible special Hermitian metric as above are given as follows (see \cite{Herm} for more details): 
\begin{enumerate}[{\rm (i)}]
\item Balanced: there exists $g_\qg$ such that $Z_{J_\qg,g_\qg}\in\zg(\kg)$ (see \eqref{koszul-g}).  

\item LCB: $\zg(\ggo)\ne 0$. i.e., $b_1(M)>0$.   

\item Pluriclosed (or SKT):  $M=G$ is a compact Lie group and $J_\ag$ is compatible with a metric $g_\ag$ such that $g_\ag|_\tg=g_b|_\tg$ for some bi-invariant metric $g_b$ on $[\ggo,\ggo]$, where $\tg$ is any maximal torus of $[\ggo,\ggo]$ (i.e., $\ag=\zg(\ggo)\oplus\tg$).  This is actually equivalent to the existence of a Bismut-flat left-invariant metric compatible with $J$.  

\item LCK: $M=G/K=G_f/K\times S^1$ (i.e., $\dim{\zg(\ggo)}=1$) and there exists $W\in\zg(\hg_f)$ such that $\ag_f=\RR W$ and the invariant ordering attached to $J_\qg$ is given by 
$$
\Delta_\qg^+:=\{\alpha\in\Delta_\qg:\alpha(\im W)>0\}. 
$$  
\end{enumerate}

The CYT condition is more involved (see \cite[Section 10]{Herm} for further details).  A necessary condition on $J=(J_\ag,J_\qg)$ for the existence of a CYT compatible metric is $Z_{J_\qg}\in\ag_f$ (i.e., $c_1(J)=0$, see \eqref{c10}), which we assume to hold from now on, and a metric $g=g_\ag+g_\qg$ is CYT if and only if 
\begin{equation}\label{CYT3}
g_\ag\left((Z_{J_\qg,g_\qg})_{\ag_f},A\right)=Q(Z_{J_\qg},A), \qquad\forall A\in\ag_f.  
\end{equation}
If $\dim{\ag_f}=1$, then this is equivalent to 
$$
g_\ag(Z_{J_\qg},Z_{J_\qg})=\tfrac{Q(Z_{J_\qg},Z_{J_\qg})^2}{Q(Z_{J_\qg,g_\qg},Z_{J_\qg})}, 
$$
so the existence always holds in that case.  

On the other hand, if in addition $J_\ag Z_{J_\qg}\in\ag_f$ (e.g., if $G$ semisimple), then another necessary condition on $J=(J_\ag,J_\qg)$ for the existence of a CYT compatible metric is  
\begin{equation}\label{CYT}
u<\tfrac{ab}{a^2+1}<v,   \qquad\mbox{for certain}\quad u<0<v \quad \mbox{depending on $J_\qg$},
\end{equation}
where $J_\ag Z_{J_\qg}=aZ_{J_\qg}+bZ_2$, $b\ne 0$, and 
$$
\ag_J:=\la Z_{J_\qg},J_\ag Z_{J_\qg}\ra_\RR =\la Z_{J_\qg},Z_2\ra_\RR, \qquad Q(Z_{J_\qg},Z_2)=0, 
\quad Q(Z_2,Z_2)=Q(Z_{J_\qg},Z_{J_\qg}).  
$$ 
Note that in terms of the basis $\{ Z_{J_\qg},Z_2\}$, 
\begin{equation}\label{Jab}
[J_\ag|_{\ag_J}] = 
\left[\begin{matrix} 
a&-\tfrac{a^2+1}{b}\\ b&-a 
\end{matrix}\right], \qquad\mbox{and}\qquad 
[g_\ag|_{\ag_J}] = 
c\left[\begin{matrix} 
b^2&-ab\\ -ab&a^2+1 
\end{matrix}\right], \quad c>0. 
\end{equation}
Condition \eqref{CYT} is also sufficient if $\ag_J=\ag_f$, e.g., when $G$ is semisimple and $\dim{\ag}=2$.

\section{Classification}\label{cla-sec}

We provide in this section a classification of C-spaces up to finite cover and equivariant diffeomorphism of dimension $\leq 10$ and a partial classification for dimension $12$.  Since the first product from the left in \eqref{dec} is always a product of complex manifolds, we will only consider C-spaces without flag factor and different from a torus.  

From a straightforward inspection of painted Dynkin diagrams (see \cite[Section 2]{Herm}), one obtains the lists of all flag manifolds of dimension $\leq 4$, $\leq 6$ and $\leq 8$, which are given in the last columns of Tables \ref{tab6}, \ref{tab8} and \ref{tab10}, respectively.

\subsection{Dimension $4$}\label{dim4-sec} 
According to \S\ref{Csp-sec}, since the only flag manifold of dimension $2$ is $S^2=\CC\PP^1=\SU(2)/\U(1)$, the only $4$-dimensional C-space up to finite cover and equivariant diffeomorphism is the Hopf manifold $M=G=\SU(2)\times S^1$.  It is easy to see using \S\ref{comp-sec} and \S\ref{herm-sec} that any left-invariant complex structure is biholomorphic to one of the form $(J(a,b),J_1)$, where $J(a,b)$ is as in \eqref{Jab}, and any Hermitian left-invariant metric descending to $F=\SU(2)/\U(1)$ is of the form $g=g_\ag+g_\qg$, where $g_\ag$ is as in \eqref{Jab} and $g_\qg=xQ|_\qg$, $x>0$.  Any complex structure $(J(a,b),J_1)$ has zero first Chern class (see \eqref{c10}) and it follows from \eqref{CYT3} that $g$ is CYT (or equivalently, Bismut flat) if and only if  $c(a^2+1)=x$.  Furthermore, any of these metrics is always pluriclosed and LCK.

\subsection{Dimension $6$}\label{dim6-sec} 
The classification of all $6$-dimensional C-spaces up to finite cover and equivariant diffeomorphism given in Table \ref{tab6} follows from the structure described in \S\ref{Csp-sec} and the list of all flag manifolds of dimension $\leq 4$.  The dimensions of the subspaces $\zg(\kg)$, $\ag_f$ and $\ag$ are provided in the table (recall that $\dim{\zg(\ggo)}=\dim{\ag}-\dim{\ag_f}$), as well as the number $r$ of irreducible isotropy summands of the flag $F$ (see \eqref{isotrep}) and the number $k$ of invariant complex structures on $F$ (see \eqref{Ccal}).  

In line 4, to obtain a closed subgroup $K$, the Lie algebra $\zg(\kg)=\kg$ of $S^1_{p,q}$ must be of the form $\RR (p\im H_\alpha+q\im H_\beta)$, where $p,q\in\ZZ$ and $\alpha,\beta$ are respectively the roots of the two copies of $\sug(2)$.  Since this is invariant under nonzero scaling, the permutation of $(p,q)$ determines a $G$-equivariant diffeomorphism and there is no flag factor, we can assume that $p,q>0$, $p\geq q$ and they are coprime.

\subsection{Dimension $8$}\label{dim8-sec}
In much the same way as for dimension $6$ in \S\ref{dim6-sec}, we obtain the classification of all $8$-dimensional C-spaces up to finite cover and equivariant diffeomorphism, which is given in Table \ref{tab8} (see also Table \ref{tab8-2}), using the list of all flag manifolds of dimension $\leq 6$.  

We give next the definitions of the isotropy subgroups $K$ in Table \ref{tab8} and some information on the flags: 
\begin{enumerate}[{\small $\bullet$}] 
\item Line 3.  The full flag $\PP(T^{1,0}\CC\PP^2)= \SU(3)/T^2$ is the projective holomorphic tangent bundle of $\CC\PP^2$.  

\item Line 5.  $S^1_{p,q}$ is as in \S\ref{dim6-sec}.  

\item Line 6.  The Lie algebra $\zg(\kg)=\kg$ of $T^2_{p,q,r}$ is the $Q$-orthogonal complement of $\ag_f=\RR (p\im H_\alpha+q\im H_\beta+r\im H_\gamma)$, where $p,q,r\in\ZZ$ and $\alpha,\beta,\gamma$ are respectively the roots of the three copies of $\sug(2)$.  Since there is no flag factor and permutations determine $G$-equivariant diffeomorphisms, we can assume that 
$$
\zg(\kg)=\la (x,y,0), (0,z,w)\ra_\RR, \qquad x,y,z,w\in\ZZ, \quad x,w\ne 0, \quad y<0, \quad z>0, 
$$
so up to the Weyl group action, we can additionally assume that $x>0$ and $w<0$.  Thus $(p,q,r)=(yw,-xw,xz)$ up to nonzero scaling, which gives $p,q,r>0$ and coprime.  Permuting again if necessary, we can finally  assume that $p\geq q\geq r$.    

\item Line 7. $\SU(2)\subset\SU(3)$ and $S^1_{p,q}$ is as in \S\ref{dim6-sec}.  

\item Line 8.  $S^1_{p,q}\subset\SU(3)$ is the Lie subgroup with Lie algebra $\RR Z_{p,q}$, where 
$$
Z_{p,q}:= p\im H_\alpha+(p+q)\im H_\beta = (p,q,-(p+q))\in\hg_f=\zg(\hg_f), \qquad p,q\in\ZZ,
$$ 
thus AW$_{p,q}:=\SU(3)/S^1_{p,q}$ is an {\it Aloff-Wallach} space.  Since this is invariant under scaling and the Weyl group permutes $(p,q,-(p+q))$ determining $\SU(3)$-equivariant diffeomorphisms, we can assume that $p>0$, $p\geq q\geq 0$ and that they are coprime.  The topology of AW$_{p,q}$ is very sensitive to the rational slope $p/q$.  Indeed, $H^4({\rm AW}_{p,q},\ZZ)=\ZZ_{p^2+pq+q^2}$ and there are infinitely many homotopy classes among them (see \cite{Krg}).  Moreover, there are different diffeomorphic classes on the same topological space (see \cite{KrcStl}).

\item Line 9.  The Stiefel manifold $V_2(\RR^5)=\Spe(2)/\SU(2)_1$ of all orthonormal $2$-frames in $\RR^5$, also presented as $\SO(5)/\SO(2)\times\SO(3)$, is simply connected and it is an $S^1$-fibration over the symmetric space  $Q^3=\Spe(2)/\U(2)$, i.e., the $3$-dimensional complex quadric.  

\item Line 10.  $S^7=\Spe(2)/\SU(2)_2$ is an $S^1$-fibration over the flag $\CC\PP^3=\Spe(2)/\U(1)\times\Spe(1)$.  

\item Line 12.  The Lie algebra of $S^1_{p,q,r}$ is $\zg(\kg)=\kg=\RR(p\im H_\alpha+q\im H_\beta+r\im H_\gamma)$, where $p,q,r\in\ZZ$ are coprime, and up to $G$-equivariant diffeomorphism, we can assume that $p,q>0$ and $p\geq q\geq r\geq 0$.
\end{enumerate}

\begin{remark}\label{B7}
We note that the Berger space $B^7=\Spe(2)/\Spe(1)_{max}$ is not an $S^1$-fibration over any flag, so $B^7\times S^1$ is not a C-space.  
\end{remark}

\subsection{Dimension $10$}\label{dim10-sec}
In Table \ref{tab10} (see also Table \ref{tab10-2}), starting from the list of all flag manifolds of dimension $\leq 8$ given in the last column and using the structure given in \S\ref{Csp-sec}, we obtain the classification, up to finite cover and equivariant diffeomorphism, of all $10$-dimensional C-spaces without flag factor which are not the product $G/K\times T^4$ or $G/K\times T^2$ for some C-space $G/K$ of dimension $6$ or $8$, respectively. 

The Lie algebras of the subgroups of $K$ denoted by $T^2_{p,q,r,s}$, $T^2_{p,q,r}$, $S^1_{p,q}$, $S^1_{p,q,r}$ and $T^2_{p,q,r,s}$ are defined in much the same way as in \S\ref{dim6-sec} and \S\ref{dim8-sec}.

\subsection{Dimension $12$}\label{dim12-sec}
We list in Table \ref{tab12} the C-spaces of dimension $12$ such that $G_f$ is simple, i.e., the flag manifold $F=G_f/H_f$, which is the base of the Tits fibration, is irreducible.  

There is a total of $21$ flags of dimension $10$: 
\begin{enumerate}[{\rm (i)}]
\item Seven irreducible ones (see Table \ref{tab12}). 

\item Four products of $\SU(3)/\U(2)$ times the flags in lines 8-11 of Table \ref{tab8}.  

\item Ten products of $\SU(2)/\U(1)$ times the flags in lines 2-9 and 11-12 of Table \ref{tab10}.  
\end{enumerate}
In order to obtain all the $12$-dimensional C-spaces which are not the product of two even dimensional C-spaces we have to consider for each of the $16$ flags of dimension $10$ with $\dim{\zg(\hg_f)}\geq 2$ and isotropy subgroup $K$ of codimension $1$ or $2$ in $H_f$, i.e., $\dim{\ag_f}=1,2$, or equivalently, $\zg(\kg)\subset\zg(\hg_f)$ of codimension $1$ or $2$.  This provides a number of $32$ C-spaces which are not the product of two even dimensional C-spaces if $\zg(\kg)$ is generically embedded in $\zg(\hg_f)$.  The only remaining C-space of dimension $12$ with this property is $\SU(2)^4/S^1_{p,q,r,s}\times S^1$.

\section{Existence of special metrics}\label{exi-sec}

We study in this section the existence of special metrics (see \S\ref{sp-sec}) on each of the invariant complex structures that admits a C-space of dimension $6$ and $8$ (see \S\ref{comp-sec}).

\subsection{Dimension $6$}\label{dim6exi-sec}
In Table \ref{tab6-exi}, the set of all $G$-invariant complex structures on each $6$-dimensional C-space in Table \ref{tab6} descending to the corresponding flag $F$ is given in the third column, which was obtained following \S\ref{comp-sec}.  The map $J_\ag(a,b):\ag\rightarrow\ag$ in line 2 is defined as in \eqref{Jab}.  The vanishing condition for the first Chern class appears in the fourth column (see \eqref{c10}).  

The rest of the columns of Table \ref{tab6-exi} indicate under what conditions there exists a special $G$-invariant Hermitian metric, which follows straightforwardly from \S\ref{sp-sec}.  An explanatory picture in the case of the complex structures in lines 4.1 and 4.2 is given in Figure \ref{6-4}.  The CYT condition on line 2 follows from \eqref{CYT} and on the rest of the lines from the fact that $\dim{\ag_f}=1$.

\subsection{Dimension $8$}\label{dim8exi-sec}
\begin{figure}
\begin{tikzpicture}[scale=.8]
\draw[->,thick, color=red] (-3.5, 0) -- (6.5, 0); 
\draw[->,thick, color=red] (0, -3.5) -- (0, 3.5); 
\draw[ultra thin,color=gray] (-3.5,-3.5) grid (6.5,3.5); 

\draw (0,2) node[right] {$\beta$};
\draw (0,2) node {{\tiny $\bullet$}};
\draw (2,0) node[below] {$\alpha$};
\draw (2,0) node {{\tiny $\bullet$}}; 
\draw (-2,0) node {{\tiny $\bullet$}};
\draw (0,-2) node {{\tiny $\bullet$}};

\draw (2,2) node[left] {$Z_{J_1}$};
\draw (2,2) node {{\tiny $\bullet$}}; 
\draw (2,-2) node[left] {$Z_{J_2}$};
\draw (2,-2) node {{\tiny $\bullet$}}; 

\draw (6,0) node[below] {$p$};
\draw (0,2) node[left] {$q$};

\draw[very thick, color=blue] (-3,-1) -- (6, 2); 
\draw (4,1) node[color=blue]{$\zg(\kg)$};

\draw[very thick, color=brown] (-1,3) -- (1,-3); 
\draw (-1,2) node[color=brown]{$\ag_f$};

\filldraw[fill=blue, opacity=0.1] (0,0) -- (6.5,0) -- (6.5,3.5) -- (3.5,3.5) -- (0,0); 
\filldraw[fill=brown, opacity=0.1] (0,0) -- (0,3.5) -- (-3.5,3.5) -- (0,0);
    
\draw (3.5,2.5) node {$J_1$};
\draw (3.5,-2.5) node {$J_2$};
\end{tikzpicture}
\caption{{\bf Dim 6}.  Proof for lines 4.1 and 4.2 in Table \ref{tab6-exi} (see \S\ref{sp-sec}).}\label{6-4}
\end{figure}
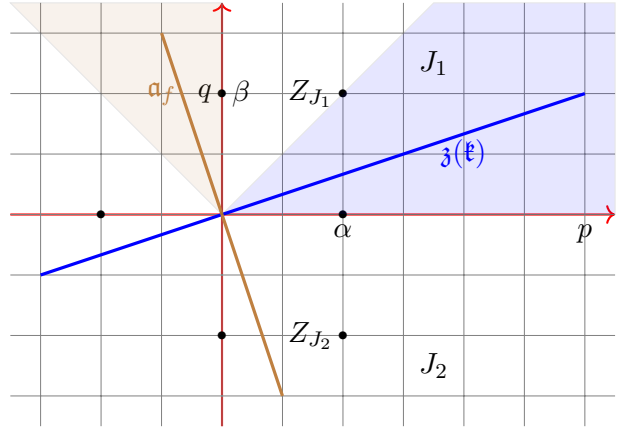

\begin{figure}
\begin{tikzpicture}[scale=.8]
\draw[->,thick] (-4.5, 0) -- (4.5, 0); 
\draw[->,thick, color=red] (0, -3.5) -- (0, 3.5); 
\draw[ultra thin,color=gray] (-4.5,-3.5) grid (4.5,3.5); 


\draw (-1,1.732) node[left] {$\beta$};
\draw (-1,1.732) node {{\tiny $\bullet$}}; 
\draw (1,-1.732) node {{\tiny $\bullet$}}; 
\draw (2,-3.464) node {{\tiny $\bullet$}};
\draw (2,-3.464) node[right] {$Z_{J_3}$};

\draw (1,1.732) node[right] {$\alpha+\beta$};
\draw (1,1.732) node {{\tiny $\bullet$}};
\draw (-1,-1.732) node {{\tiny $\bullet$}};
\draw (2,3.464) node {{\tiny $\bullet$}};
\draw (2,3.464) node[left] {$Z_{J_1}$};

\draw (2,0) node[below] {$\alpha$};
\draw (2,0) node {{\tiny $\bullet$}}; 
\draw (-2,0) node {{\tiny $\bullet$}}; 
\draw (4,0) node {{\tiny $\bullet$}};
\draw (4,0) node[below] {$Z_{J_2}$};


\draw[thick, color=red] (-4,2.464) -- (4,-2.464); 
\draw[thick, color=red] (-4,-2.464) -- (4,2.464); 

\draw[very thick, color=brown] (-4,-1.33) -- (4, 1.33); 
\draw (3.5,1.5) node[color=brown]{$\ag_f$};
\filldraw[fill=brown, opacity=0.1] (0,0) -- (4,2.464) -- (4,-2.464);

\draw[very thick, color=blue] (-1,3) -- (1,-3); 
\draw (-0.35,2.6) node[color=blue]{$\zg(\kg)$};
\draw[thick, color=blue, opacity=0.25] (2,3.464) -- (-2,-3.464); 
\filldraw[fill=blue, opacity=0.1] (0,0) -- (2,3.464) -- (-2,3.464);

\draw (2,-3.464) -- (-2,3.464); 

\draw (2.5,2.5) node {$J_1$};
\draw (3.5,0.5) node {$J_2$};
\draw (2.5,-2.5) node {$J_3$};
\end{tikzpicture}
\caption{{\bf Dim 8}. Proof for lines 8.1-8.3 in Table \ref{tab8-exi} (see \S\ref{sp-sec}).}\label{8-8}
\end{figure}
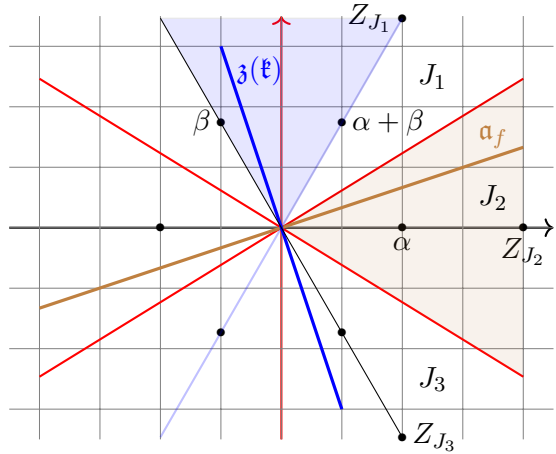

The $G$-invariant complex structures on each C-space of dimension $8$ in Table \ref{tab8} descending to the corresponding flag $F$ are given in Table \ref{tab8-exi} (see \S\ref{comp-sec}).  In lines 3, 12.1-12.4 and 13, the map $J_\ag(a,b)$ is as in \eqref{Jab}, and in lines 6.1-6.4 and 12.1-12.4 the invariant orderings defining $J_1$, $J_2$, $J_3$ and $J_4$ are respectively generated by the following sets of simple roots: 
$$
\{\alpha,\beta,\gamma\}, \quad \{\alpha,-\beta,\gamma\}, \quad \{-\alpha,-\beta,\gamma\}, \quad \{-\alpha,\beta,\gamma\}.  
$$
The corresponding Koszul vectors are therefore given by 
$$
Z_{J_1}=(1,1,1), \quad Z_{J_2}=(1,-1,1), \quad Z_{J_3}=(-1,-1,1), \quad Z_{J_4}=(-1,1,1).
$$
For lines 5.1-5.2 and 7.1-7.2, the $J_i$ are defined as in Figure \ref{6-4} and for lines 8.1-8.3 as in Figure \ref{8-8}.  
 
In Table \ref{tab8-exi}, it is also described under what conditions on the C-space and $J$ there exists a special $G$-invariant Hermitian metric (see \S\ref{sp-sec}).  An explanatory picture in the case of the space in lines 8.1-8.3 is given in Figure \ref{8-8}. 

The CYT condition on all the lines excepting 2, 3, 12 and 13 follows from the fact that $\dim{\ag_f}=1$ (see Table \ref{tab8-2}).  On lines 3, 12.1-12.4 and 13, the CYT condition follows from \eqref{CYT}.  For the Lie group $M=\SU(2)^2\times T^2$ in line 2, we use the fact that either $J_\ag\ag_f=\ag_f$ and $\ag_J=\ag_f$, or $J_\ag\ag_f\ne\ag_f$ and hence $J_\ag\ag_f\cap\ag_f=0$, which implies that $\ag=\ag_f\oplus J_\ag\ag_f$ and so any metric $g$ such that $q_\qg=Q|_\qg$ and $\{\im H_\alpha, \im H_\beta, J_\ag\im H_\alpha, J_\ag\im H_\beta\}$ is $g_\ag$ orthonormal is CYT by \eqref{CYT3}.

\begin{table}
\centering
\begin{tabular}{c|c | c | c | c | c | c | c | cc}
\hline\hline \rule{0pt}{14pt}
& $G/K$ & $M$ & $\zg(\kg)$ & $\ag_f$ & $\ag$ & $r$ & $k$ & $F=G_f/H_f$ & \\ [2pt]
\hline\hline\rule{0pt}{14pt}
1& $\SU(2) \times T^{3}$ & $S^{3} \times T^{3}$ & $0$ & $1$ & $4$ & $1$ & $2$ & $S^2=\CC\PP^1=\tfrac{\SU(2)}{\U(1)}$ & \dynkin[scale=1.5]{A}{o} \\ 
\hline \rule{0pt}{14pt}
2& $\SU(2)^2$ & $S^{3} \times S^{3}$ & $0$ & $2$ & $2$ & $2$ & $4$ & $S^2\times S^2$ & \dynkin[scale=1.5] {A}{o} \dynkin[scale=1.5]{A}{o} \\ 
\hline \rule{0pt}{14pt}
3& $\frac{\SU(3)}{\SU(2)} \times S^{1}$ & $S^{5} \times S^{1}$ & $0$ & $1$ & $2$ & $1$ & $2$ & $\CC\PP^2=\frac{\SU(3)} {\U(2)}$ & \dynkin[scale=1.5]{A}{o*} \\ 
\hline \rule{0pt}{14pt}
4& $\frac{\SU(2)^2}{S^1_{p,q}} \times S^{1}$ & $S^{3} \times S^{2} \times S^{1}$ & $1$ & $1$ & $2$ & $2$ & $4$ & $S^2\times S^2$ & \dynkin[scale=1.5]{A}{o} \dynkin[scale=1.5]{A}{o} \\ 
\hline\hline 
\end{tabular}
\vspace{.3cm}
\caption{{\bf Dim 6}.  C-spaces of dimension $6$ without flag factor.  In line 4, $p,q\in\NN$, $p\geq q$ and  they are coprime.}
\label{tab6}
\end{table}

\begin{table}
\centering
\begin{tabular}{c|c | c | c | cc}
\hline\hline \rule{0pt}{14pt}
& $G/K$ & $M$ & $\dim{\mathfrak{a}}$ & $F=G_f/H_f$ & \\ [2pt]
\hline\hline 
\rule{0pt}{14pt}
1& $\SU(2) \times T^{5}$ & $S^{3} \times T^{5}$ & $6$ & $S^2$ & \dynkin[scale=1.5]{A}{o} \\ 
\hline \rule{0pt}{14pt}
2& $\SU(2)^2 \times T^{2}$ & $S^{3} \times S^{3} \times T^{2}$ & $4$ & $S^2\times S^2$ & \dynkin[scale=1.5]{A}{o} \dynkin[scale=1.5]{A}{o} \\ 
\hline \rule{0pt}{14pt}
3& $\SU(3)$ & -- & $2$ & $\tfrac{\SU(3)}{T^2}$ & \dynkin[scale=1.5]{A}{oo} \\ 
\hline \rule{0pt}{14pt}
4& $\frac{\SU(3)}{\SU(2)} \times T^{3}$ & $S^{5} \times T^{3}$ & $4$ & $\CC\PP^2=\frac{\SU(3)} {\U(2)}$ & \dynkin[scale=1.5]{A}{o*} \\ 
\hline \rule{0pt}{14pt}
5& $\frac{\SU(2)^2}{S^1_{p,q}} \times T^{3}$, $p,q\in\NN$ & $S^{3} \times S^{2} \times T^{3}$ & $4$ & $S^2\times S^2$ & \dynkin[scale=1.5]{A}{o} \dynkin[scale=1.5]{A}{o} \\ 
\hline \rule{0pt}{14pt}
6& $\frac{\SU(2)^3}{T^2_{p,q,r}} \times S^{1}$, $p,q,r\in\NN$ & $S^{3} \times S^{2} \times S^{2}\times S^1$ & $2$ & $S^2\times S^2\times S^2$ & \dynkin[scale=1.5]{A}{o} \dynkin[scale=1.5]{A}{o} \dynkin[scale=1.5]{A}{o} \\ 
\hline \rule{0pt}{14pt}
7& $\frac{\SU(3)\times \SU(2)}{\SU(2)\times S^1_{p,q}} \times S^1$ & $S^5\times S^{2}\times S^1$ & $2$ & $\CC\PP^2\times S^2$ & \dynkin[scale=1.5]{A}{o*} \dynkin[scale=1.5]{A}{o} \\ 
\hline \rule{0pt}{14pt}
8& $\frac{\SU(3)}{S^1_{p,q}} \times S^1$, $p\in\NN$, $q\in\NN_0$ & AW$_{p,q}$ $\times S^1$ & $2$ & $\tfrac{\SU(3)}{T^2}$ & \dynkin[scale=1.5]{A}{oo} \\ 
\hline \rule{0pt}{14pt}
9& $\frac{\Spe(2)}{\SU(2)_1} \times S^{1}$ & $V_2(\RR^5)\times S^1$ & $2$ & $Q^3=\tfrac{\Spe(2)}{\U(2)}$ & \dynkin[scale=1.5]{B}{o*} \\ 
\hline \rule{0pt}{14pt}
10& $\frac{\Spe(2)}{\SU(2)_2} \times S^{1}$ & $S^7\times S^1$ & $2$ & $\CC\PP^3=\tfrac{\Spe(2)}{\U(1)\times\Spe(1)}$ & \dynkin[scale=1.5]{B}{*o} \\ 
\hline \rule{0pt}{14pt}
11& $\frac{\SU(4)}{\SU(3)} \times S^{1}$ & $S^{7} \times S^{1}$ & $2$ & $\CC\PP^3=\frac{\SU(4)} {\U(3)}$ & \dynkin[scale=1.5]{A}{o**} \\ 
\hline \rule{0pt}{14pt}
12& $\frac{\SU(2)^3}{S^1_{p,q,r}}$, $p,q\in\NN$, $r\in\NN_0$ & $S^{3} \times S^3 \times S^{2}$ & $2$ & $S^2\times S^2\times S^2$ & \dynkin[scale=1.5]{A}{o} \dynkin[scale=1.5]{A}{o} \dynkin[scale=1.5]{A}{o} \\ 
\hline \rule{0pt}{14pt}
13& $\frac{\SU(3)}{\SU(2)} \times \SU(2)$ & $S^{5} \times S^{3}$ & $2$ & $\CC\PP^2\times S^2$ & \dynkin[scale=1.5]{A}{o*} \dynkin[scale=1.5]{A}{o} \\ 
\hline\hline
\end{tabular}
\vspace{.3cm}
\caption{{\bf Dim 8}.  C-spaces of dimension $8$ without flag factor (see \S\ref{dim8-sec} for extra conditions on the integers $p,q,r$).}
\label{tab8}
\end{table}

\begin{table}
\centering
\begin{tabular}{c|c | c | c | c | c | c}
\hline\hline \rule{0pt}{14pt}
& $G/K$ & $\dim{\zg(\mathfrak{k})}$ & $\dim{\mathfrak{a}_f}$ & $\dim{\mathfrak{a}}$ & $r$ & $k$ \\ [2pt]
\hline\hline 
\rule{0pt}{14pt}
1& $\SU(2) \times T^{5}$ & $0$ & $1$ & $6$ & $1$ & $2$ \\ 
\hline \rule{0pt}{14pt}
2& $\SU(2)^2 \times T^{2}$ & $0$ & $2$ & $4$ & $2$ & $4$ \\ 
\hline \rule{0pt}{14pt}
3& $\SU(3)$ & $0$ & $2$ & $2$ & $3$ & $6$  \\ 
\hline \rule{0pt}{14pt}
4& $\frac{\SU(3)}{\SU(2)} \times T^{3}$ & $0$ & $1$ & $4$ & $1$ & $2$ \\ 
\hline \rule{0pt}{14pt}
5& $\frac{\SU(2)^2}{S^1_{p,q}} \times T^{3}$ & $1$ & $1$ & $4$ & $2$ & $4$ \\ 
\hline \rule{0pt}{14pt}
6& $\frac{\SU(2)^3}{T^2_{p,q,r}} \times S^{1}$ & $2$ & $1$ & $2$ & $3$ & $8$ \\ 
\hline \rule{0pt}{14pt}
7& $\frac{\SU(3)\times \SU(2)}{\SU(2)\times S^1_{p,q}} \times S^1$ & $1$ & $1$ & $2$ & $2$ & $4$ \\ 
\hline \rule{0pt}{14pt}
8& $\frac{\SU(3)}{S^1_{p,q}} \times S^1$ & $1$ & $1$ & $2$ & $3$ & $6$ \\ 
\hline \rule{0pt}{14pt}
9& $\frac{\Spe(2)}{\SU(2)_1} \times S^{1}$ & $0$ & $1$ & $2$ & $1$ & $2$ \\ 
\hline \rule{0pt}{14pt}
10& $\frac{\Spe(2)}{\SU(2)_2} \times S^{1}$ & $0$ & $1$ & $2$ & $2$ & $2$ \\ 
\hline \rule{0pt}{14pt}
11& $\frac{\SU(4)}{\SU(3)} \times S^{1}$ & $0$ & $1$ & $2$ & $1$ & $2$ \\ 
\hline \rule{0pt}{14pt}
12& $\frac{\SU(2)^3}{S^1_{p,q,r}}$ & $1$ & $2$ & $2$ & $3$ & $8$ \\ 
\hline \rule{0pt}{14pt}
13& $\frac{\SU(3)}{\SU(2)} \times \SU(2)$ & $0$ & $2$ & $2$ & $2$ & $4$ \\ 
\hline\hline
\end{tabular}
\vspace{.3cm}
\caption{{\bf Dim 8}.  C-spaces of dimension $8$ without flag factor.}
\label{tab8-2}
\end{table}

\begin{table}
\centering
\begin{tabular}{c|c | c | c | c}
\hline\hline \rule{0pt}{14pt}
& $G/K$ & $M$ & $\dim{\mathfrak{a}}$ & $F$ \\ [2pt]
\hline\hline \rule{0pt}{14pt}
1& $\SU(2)^{3} \times S^{1}$ & $S^3 \times S^3 \times S^3 \times S^1$ & $4$ & \dynkin[scale=1.5]{A}{o} \dynkin[scale=1.5]{A}{o} \dynkin[scale=1.5]{A}{o} \\ 
\hline \rule{0pt}{14pt}
2& $\Spe(2)$ & -- & $2$ & $\tfrac{\Spe(2)}{T^2}$ \dynkin[scale=1.5]{B}{oo} \\ 
\hline \rule{0pt}{14pt}
3& $\frac{\SU(2)^{4}}{T^3_{p,q,r,s}} \times S^{1}$ & $S^3\times S^2\times S^2\times S^2\times S^1$ & $2$ & \dynkin[scale=1.5]{A}{o} \dynkin[scale=1.5]{A}{o} \dynkin[scale=1.5]{A}{o} \dynkin[scale=1.5]{A}{o} \\ 
\hline \rule{0pt}{14pt}
4& $\frac{\SU(3) \times \SU(2)^2}{\SU(2) \times T_{p,q,r}^2} \times S^{1}$ & $S^5\times S^2\times S^2\times S^1$ & $2$ & \dynkin[scale=1.5]{A}{o*} \dynkin[scale=1.5]{A}{o} \dynkin[scale=1.5]{A}{o} \\ 
\hline \rule{0pt}{14pt}
5& $\frac{\SU(3) \times \SU(2)}{T_{p,q,r}^2} \times S^{1}$ & AW$_{p,q}\times S^2\times S^1$ & $2$ & \dynkin[scale=1.5]{A}{oo} \dynkin[scale=1.5]{A}{o} \\ 
\hline \rule{0pt}{14pt}
6& $\tfrac{\Spe(2) \times \SU(2)}{\SU(2)_1 \times S_{p,q}^1} \times S^{1}$ & $V_2(\RR^5)\times S^2\times S^1$ & $2$ & \dynkin[scale=1.5]{B}{o*} \dynkin[scale=1.5]{A}{o} \\ 
\hline \rule{0pt}{14pt}
7& $\textstyle\frac{\Spe(2) \times \SU(2)}{\SU(2)_2 \times S_{p,q}^1} \times S^{1}$ & $S^7\times S^2\times S^1$ & $2$ & \dynkin[scale=1.5]{B}{*o} \dynkin[scale=1.5]{A}{o} \\ 
\hline \rule{0pt}{14pt}
8& $\textstyle\frac{\SU(4) \times \SU(2)}{\SU(3) \times S_{p,q}^1} \times S^{1}$ & $S^7\times S^2\times S^1$ & $2$ & \dynkin[scale=1.5]{A}{o**} \dynkin[scale=1.5]{A}{o} \\ 
\hline \rule{0pt}{14pt}
9& $\tfrac{\SU(3) \times \SU(3)}{\SU(2) \times \SU(2) \times S_{p,q}^1} \times S^{1}$ & -- & $2$ & \dynkin[scale=1.5]{A}{o*} \dynkin[scale=1.5]{A}{o*} \\ 
\hline \rule{0pt}{14pt}
10& $\tfrac{\Spe(2)}{S_{p,q}^1} \times S^{1}$ & -- & $2$ & \dynkin[scale=1.5]{B}{oo} \\ 
\hline \rule{0pt}{14pt}
11& $\tfrac{\SU(4)}{\SU(2) \times \SU(2)} \times S^{1}$ & -- & $2$ & \dynkin[scale=1.5]{A}{*o*} \\ 
\hline \rule{0pt}{14pt}
12& $\frac{\SU(5)}{\SU(4)} \times S^{1}$ & $S^{9} \times S^{1}$ & $2$ & $\CC\PP^4$ \dynkin[scale=1.5]{A}{o***} \\ 
\hline \rule{0pt}{14pt}
13& $\frac{\SU(3)}{\SU(2)}\times \tfrac{\SU(3)}{\SU(2)}$ & $S^5 \times S^5$ & $2$ & \dynkin[scale=1.5]{A}{o*} \dynkin[scale=1.5]{A}{o*} \\ 
\hline \rule{0pt}{14pt}
14& $\frac{\SU(3) \times \SU(2)^2}{\SU(2) \times S^1_{p,q,r}}$ & $S^5 \times S^3 \times S^2$ & $2$ & \dynkin[scale=1.5]{A}{o*} \dynkin[scale=1.5]{A}{o} \dynkin[scale=1.5]{A}{o} \\ 
\hline \rule{0pt}{14pt}
15& $\frac{\Spe(2)}{\SU(2)_1} \times \SU(2) $ & $V_2(\RR^5)\times S^3$ & $2$ & \dynkin[scale=1.5]{B}{o*} \dynkin[scale=1.5]{A}{o} \\ 
\hline \rule{0pt}{14pt}
16& $\frac{\Spe(2)}{\SU(2)_2} \times \SU(2) $ & $S^7\times S^3$ & $2$ & \dynkin[scale=1.5]{B}{*o} \dynkin[scale=1.5]{A}{o} \\ 
\hline \rule{0pt}{14pt}
17& $\frac{\SU(4)}{\SU(3)} \times \SU(2) $ & $S^7 \times S^3$ & $2$ & \dynkin[scale=1.5]{A}{o**} \dynkin[scale=1.5]{A}{o} \\ 
\hline \rule{0pt}{14pt}
18& $\tfrac{\SU(2)^4}{T^2_{p,q,r,s}} $ & $S^3\times S^3\times S^2\times S^2$ & $2$ & \dynkin[scale=1.5]{A}{o} \dynkin[scale=1.5]{A}{o} \dynkin[scale=1.5]{A}{o} \dynkin[scale=1.5]{A}{o} \\ 
\hline \rule{0pt}{14pt}
19& $\tfrac{\SU(3) \times \SU(2)}{S_{p,q,r}^1}$ & AW$_{p,q}\times S^3$ & $2$ & \dynkin[scale=1.5]{A}{oo} \dynkin[scale=1.5]{A}{o} \\ 
\hline\hline
\end{tabular}
\vspace{.3cm}
\caption{{\bf Dim 10}.  C-spaces of dimension $10$ without flag factor which are not the product $G/K\times T^4$ or $G/K\times T^2$ for some C-space $G/K$ of dimension $6$ or $8$, respectively. }
\label{tab10}
\end{table}

\begin{table}
\centering
\begin{tabular}{c|c | c | c | c | c | c}
\hline\hline \rule{0pt}{14pt}
& $G/K$ & $\dim{\zg(\mathfrak{k})}$ & $\dim{\mathfrak{a}_f}$ & $\dim{\mathfrak{a}}$ & $r$ & $k$ \\ [2pt]
\hline\hline \rule{0pt}{14pt}
1& $\SU(2)^{3} \times S^{1}$ & $0$ & $3$ & $4$ & $3$ & $8$ \\ 
\hline \rule{0pt}{14pt}
2& $\Spe(2)$ & $0$ & $2$ & $2$ & $4$ & $8$ \\ 
\hline \rule{0pt}{14pt}
3& $\frac{\SU(2)^{4}}{T^3_{p,q,r,s}} \times S^{1}$ & $3$ & $1$ & $2$ & $4$ & $16$ \\ 
\hline \rule{0pt}{14pt}
4& $\frac{\SU(3) \times \SU(2)^2}{\SU(2) \times T_{p,q,r}^2} \times S^{1}$ & $2$ & $1$ & $2$ & $3$ & $8$ \\ 
\hline \rule{0pt}{14pt}
5& $\frac{\SU(3) \times \SU(2)}{T_{p,q,r}^2} \times S^{1}$ & $2$ & $1$ & $2$ & $4$ & $12$ \\ 
\hline \rule{0pt}{14pt}
6& $\textstyle\frac{\Spe(2) \times \SU(2)}{\SU(2)_1 \times S_{p,q}^1} \times S^{1}$ & $1$ & $1$ & $2$ & $2$ & $4$ \\ 
\hline \rule{0pt}{14pt}
7& $\textstyle\frac{\Spe(2) \times \SU(2)}{\SU(2)_2 \times S_{p,q}^1} \times S^{1}$ & $1$ & $1$ & $2$ & $3$ & $4$ \\ 
\hline \rule{0pt}{14pt}
8& $\textstyle\frac{\SU(4) \times \SU(2)}{\SU(3) \times S_{p,q}^1} \times S^{1}$ & $1$ & $1$ & $2$ & $2$ & $4$ \\ 
\hline \rule{0pt}{14pt}
9& $\textstyle\frac{\SU(3) \times \SU(3)}{\SU(2) \times \SU(2) \times S_{p,q}^1} \times S^{1}$ & $1$ & $1$ & $2$ & $2$ & $4$ \\ 
\hline \rule{0pt}{14pt}
10& $\textstyle\frac{\Spe(2)}{S_{p,q}^1} \times S^{1}$ & $1$ & $1$ & $2$ & $4$ & $8$ \\ 
\hline \rule{0pt}{14pt}
11& $\textstyle\frac{\SU(4)}{\SU(2) \times \SU(2)} \times S^{1}$ & $0$ & $1$ & $2$ & $1$ & $2$ \\ 
\hline \rule{0pt}{14pt}
12& $\frac{\SU(5)}{\SU(4)} \times S^{1}$ & $0$ & $1$ & $2$ & $1$ & $2$ \\ 
\hline \rule{0pt}{14pt}
13& $\frac{\SU(3)}{\SU(2)} \times \frac{\SU(3)}{\SU(2)}$ & $0$ & $2$ & $2$ & $2$ & $4$ \\ 
\hline \rule{0pt}{14pt}
14& $\frac{\SU(3) \times \SU(2)^2}{\SU(2) \times S^1_{p,q,r}}$ & $1$ & $2$ & $2$ & $3$ & $8$ \\ 
\hline \rule{0pt}{14pt}
15& $\frac{\Spe(2)}{\SU(2)_1} \times \SU(2) $ & $0$ & $2$ & $2$ & $2$ & $4$ \\ 
\hline \rule{0pt}{14pt}
16& $\frac{\Spe(2)}{\SU(2)_2} \times \SU(2) $ & $0$ & $2$ & $2$ & $3$ & $4$ \\ 
\hline \rule{0pt}{14pt}
17& $\frac{\SU(4)}{\SU(3)} \times \SU(2) $ & $0$ & $2$ & $2$ & $2$ & $4$ \\ 
\hline \rule{0pt}{14pt}
18& $\textstyle\frac{\SU(2)^4}{T^2_{p,q,r,s}} $ & $2$ &$2$ & $2$ & $4$ & $16$ \\ 
\hline \rule{0pt}{14pt}
19& $\textstyle\frac{\SU(3) \times \SU(2)}{S_{p,q,r}^1}$ & $1$ & $2$ & $2$ & $4$ & $12$ \\ 
\hline\hline
\end{tabular}
\vspace{.3cm}
\caption{{\bf Dim 10}.  C-spaces of dimension $10$ without flag factor which are not the product $G/K\times T^4$ or $G/K\times T^2$ for some C-space $G/K$ of dimension $6$ or $8$, respectively.}
\label{tab10-2}
\end{table}

\begin{table}
\centering
\begin{tabular}{c | c | c | c | c | cc}
\hline\hline \rule{0pt}{14pt}
& $G/K$ & $\ag_f$ & $r$ & $k$ & $F=G_f/H_f$ & \\ [2pt] \hline
\hline\rule{0pt}{14pt}
1& $\tfrac{\SU(6)}{\SU(5)}\times S^1$ & $1$ & $1$ & $2$ & $\CC\PP^5=\tfrac{\SU(6)}{\U(5)}$ & \dynkin[scale=1.5]{A}{o*****} \\ 
\hline \rule{0pt}{14pt}
2& $\tfrac{\SU(4)}{\SU(2)_1}$ & $2$ & $3$ & $6$ & $\tfrac{\SU(4)}{\Se(\U(1)\times\U(1)\times\U(2))}$ & \dynkin[scale=1.5] {A}{oo*} \\ 
\hline \rule{0pt}{14pt}
3& $\tfrac{\SU(4)}{\SU(2)_1\times S^1_{p,q}}\times S^1$ & $1$ & $3$ & $6$ & $\tfrac{\SU(4)}{\Se(\U(1)\times\U(1)\times\U(2))}$ & \dynkin[scale=1.5] {A}{oo*} \\ 
%
\hline \rule{0pt}{14pt}
4& $\tfrac{\SU(4)}{\SU(2)_2\times S^1_{p,q}}\times S^1$ & $1$ & $3$ & $6$ & $\tfrac{\SU(4)}{\Se(\U(1)\times\U(2)\times\U(1))}$ & \dynkin[scale=1.5] {A}{o*o} \\ 
\hline \rule{0pt}{14pt}
5& $\tfrac{\SO(7)}{\SO(5)} \times S^{1}$ & $1$ & $1$ & $2$ & $V_2(\RR^7)=\tfrac{\SO(7)}{\SO(2)\times\SO(5)}$ & \dynkin[scale=1.5]{B}{o**}  \\ 
\hline \rule{0pt}{14pt}
6& $\tfrac{\Spe(3)}{\Spe(2)} \times S^{1}$ & $1$ & $2$ & $2$ & $\tfrac{\Spe(3)}{S^1\times\Spe(2)}$ & \dynkin[scale=1.5]{C}{o**}  \\ 
\hline \rule{0pt}{14pt}
7& $\tfrac{\Gg_2}{\SU(2)_1} \times S^{1}$ & $1$ & $2$ & $2$ & $\tfrac{\Gg_2}{\U(2)_1}$ & \dynkin[scale=1.5]{G}{o*}  \\ 
\hline \rule{0pt}{14pt}
8& $\tfrac{\Gg_2}{\SU(2)_2} \times S^{1}$ & $1$ & $3$ & $2$ & $\tfrac{\Gg_2}{\U(2)_2}$ & \dynkin[scale=1.5]{G}{*o}  \\ 
\hline\hline 
\end{tabular}
\vspace{.3cm}
\caption{{\bf Dim 12}.  C-spaces of dimension $12$ such that $G_f$ is simple (i.e., the flag $F=G_f/H_f$ is irreducible).  In all the lines $\dim{\ag}=2$. }
\label{tab12}
\end{table}

\begin{table}
\centering
\begin{tabular}{c|c | c | c | c | c | c | c | c}
\hline\hline \rule{0pt}{14pt}
&$G/K$ & $J$ & $c_1(J) = 0$ & Bal & LCB & SKT & CYT & LCK \\ [2pt]
\hline\hline \rule{0pt}{14pt}
1& $\SU(2) \times T^{3}$ & $(J_\ag,J_i)$ & yes & no & yes & yes & yes & no \\ 
\hline \rule{0pt}{14pt}
2& $\SU(2)^2$ & $(J_\ag(a,b),J_i)$ & yes & no & no & $a = 0$ & $\tfrac{|ab|}{a^2+1} < 1$ & no \\ 
\hline \rule{0pt}{14pt}
3& $\frac{\SU(3)}{\SU(2)} \times S^{1}$ & $(J_\ag,J_i)$ & yes & no & yes & no & yes & yes \\ 
\hline \rule{0pt}{14pt}
4.1& $\frac{\SU(2)^2}{S^1_{p,q}} \times S^{1}$ & $(J_\ag,\pm J_1)$ & no & yes & yes & no & no & no \\ 
\hline \rule{0pt}{14pt}
4.2& $\frac{\SU(2)^2}{S^1_{p,q}} \times S^{1}$ & $(J_\ag,\pm J_2)$ & $p = q=1$ & no & yes & no & $p = q=1$ & yes \\ 
\hline\hline
\end{tabular}
\vspace{.3cm}
\caption{{\bf Dim 6}.  Complex structures on $6$-dimensional C-spaces as in Table \ref{tab6} and the existence conditions for special metrics.}
\label{tab6-exi}
\end{table}

\begin{table}
\centering
\begin{tabular}{c|c | c | c | c | c | c | c}
\hline\hline \rule{0pt}{14pt}
& $(G/K,J)$ & $J$ & $c_1(J) = 0$ & Bal & SKT & CYT & LCK \\ [2pt]
\hline\hline \rule{0pt}{14pt}
1& $\SU(2) \times T^{5}$ & $(J_\ag,J_i)$ & yes & no & yes & yes & no \\ 
\hline \rule{0pt}{14pt}
2.1& $\SU(2)^2 \times T^{2}$ & $\begin{matrix}(J_\ag,J_i),\\ J_\ag\ag_f=\ag_f \end{matrix}$ & yes & no & (1) & (2) & no \\ 
\hline \rule{0pt}{14pt}
2.2& $\SU(2)^2 \times T^{2}$ & $\begin{matrix}(J_\ag,J_i),\\ J_\ag\ag_f\ne \ag_f \end{matrix}$ & yes & no & (1) & yes & no \\ 
\hline \rule{0pt}{14pt}
3& $\SU(3)$ & $(J_\ag(a,b),J_i)$ & yes & no & (3) & (2) & no \\ 
\hline \rule{0pt}{14pt}
4& $\frac{\SU(3)}{\SU(2)} \times T^{3}$ & $(J_\ag,J_i)$ & yes & no & no & yes & no \\ 
\hline \rule{0pt}{14pt}
5.1& $\frac{\SU(2)^2}{S^1_{p,q}} \times T^{3}$ & $(J_\ag,\pm J_1)$ & no & yes & no & no & no \\ 
\hline \rule{0pt}{14pt}
5.2& $\frac{\SU(2)^2}{S^1_{p,q}} \times T^{3}$ & $(J_\ag,\pm J_2)$ & $p=q=1$ & no & no & $p=q=1$ & no \\ 
\hline \rule{0pt}{14pt}
6.1& $\frac{\SU(2)^{3}}{T^2_{p,q,r}} \times S^{1}$ & $(J_\ag,\pm J_1)$ & $\begin{matrix} p=q=1,\\ r=1 \end{matrix}$ & no & no & $\begin{matrix} p=q=1,\\ r=1 \end{matrix}$ & yes \\ 
\hline \rule{0pt}{14pt}
6.2& $\frac{\SU(2)^{3}}{T^2_{p,q,r}} \times S^{1}$ & $(J_\ag,\pm J_2)$ & no & yes & no & no & no \\ 
\hline \rule{0pt}{14pt}
6.3& $\frac{\SU(2)^{3}}{T^2_{p,q,r}} \times S^{1}$ & $(J_\ag,\pm J_3)$ & no & yes & no & no & no \\ 
\hline \rule{0pt}{14pt}
6.4& $\frac{\SU(2)^{3}}{T^2_{p,q,r}} \times S^{1}$ & $(J_\ag,\pm J_4)$ & no & yes & no & no & no \\ 
\hline \rule{0pt}{14pt}
7.1& $\frac{\SU(3)\times \SU(2)}{\SU(2)\times S^1_{p,q}} \times S^1$ & $(J_\ag,\pm J_1)$ & no & yes & no & no & no \\ 
\hline \rule{0pt}{14pt}
7.2& $\frac{\SU(3)\times \SU(2)}{\SU(2)\times S^1_{p,q}} \times S^1$ & $(J_\ag,\pm J_2)$ & $p = q=1$ & no & no & $p = q=1$ & yes \\ 
\hline \rule{0pt}{14pt}
8.1& $\frac{\SU(3)}{S^1_{p,q}} \times S^1$ & $(J_\ag,\pm J_1)$ & no & yes & no & no & no \\ 
\hline \rule{0pt}{14pt}
8.2& $\frac{\SU(3)}{S^1_{p,q}} \times S^1$ & $(J_\ag,\pm J_2)$ & $p=q=1$ & no & no & $p=q=1$ & $q\in\NN$ \\ 
\hline \rule{0pt}{14pt}
8.3& $\frac{\SU(3)}{S^1_{p,q}} \times S^1$ & $(J_\ag,\pm J_3)$ & no & yes & no & no & no \\ 
\hline \rule{0pt}{14pt}
9& $\frac{\Spe(2)}{\SU(2)_1} \times S^{1}$ & $(J_\ag,J_i)$ & yes & no & no & yes & yes \\ 
\hline \rule{0pt}{14pt}
10& $\frac{\Spe(2)}{\SU(2)_2} \times S^{1}$ & $(J_\ag,J_i)$ & yes & no & no & yes & yes \\ 
\hline \rule{0pt}{14pt}
11& $\frac{\SU(4)}{\SU(3)} \times S^{1}$ & $(J_\ag,J_i)$ & yes & no & no & yes & yes \\ 
\hline \rule{0pt}{14pt}
12.1& $\frac{\SU(2)^{3}}{S^1_{p,q,r}}$ & $(J_\ag(a,b),\pm J_1)$ & no & $r\in\NN$ & no & no & no \\ 
\hline \rule{0pt}{14pt}
12.2& $\frac{\SU(2)^{3}}{S^1_{p,q,r}}$ & $(J_\ag(a,b),\pm J_2)$ & $\begin{matrix} p=q=1,\\ r=0 \end{matrix}$ & no & no & (2)  & no \\ 
\hline \rule{0pt}{14pt}
12.3& $\frac{\SU(2)^{3}}{S^1_{p,q,r}}$ & $(J_\ag(a,b),\pm J_3)$ & no & no & no & no & no \\ 
\hline \rule{0pt}{14pt}
12.4& $\frac{\SU(2)^{3}}{S^1_{p,q,r}}$ & $(J_\ag(a,b),\pm J_4)$ & $p=q+r$  & no & no & (2) & no \\ 
\hline \rule{0pt}{14pt}
13& $\frac{\SU(3)}{\SU(2)} \times \SU(2)$ & $(J_\ag(a,b),\pm J_i)$ & yes & no & no & (2) & no \\ 
\hline\hline
\end{tabular}
\vspace{.3cm}
\caption{{\bf Dim 8}.  Existence conditions for special metrics on $8$-dimensional C-spaces without flag factor. (1): $J_\ag$ is compatible with the restriction on $\ag$ of a bi-invariant metric of $G$.  (2): $u < \frac{ab}{a^2+1} < v$ for suitable $u<0<v$ depending on $J_i$ and the condition for $c_1(J)=0$ holds.  (3): $a = 0$, $b = \pm 1$.}
\label{tab8-exi}
\end{table}

\end{document}